
\magnification=1100
\baselineskip=15pt

\def\ref#1{{\rm [}{\bf #1}{\rm ]}}   

\outer\def\proclaim#1{\medbreak\noindent\bf\ignorespaces
   #1\unskip.\enspace\sl\ignorespaces}
\outer\def\endproclaim{\par\ifdim\lastskip<\medskipamount\removelastskip
   \penalty 55 \fi\medskip\rm}

\def\EEE{{\bf E}}       
      
      \def\PPP{{\bf P}}

\def\II{{\cal I}}

\def\nin{\noindent}

\def\rect#1#2#3{\raise .1ex\vbox{\hrule height.#3pt
   \hbox{\vrule width.#3pt height#2pt \kern#1pt\vrule width.#3pt}
        \hrule height.#3pt}}

\def\qed{$\hskip 5pt\rect364$} 

\centerline{\bf On a Proof of the Strong Law of Large Numbers}
\medskip
\centerline{P. J. Fitzsimmons}
\centerline{Department of Mathematics}
\centerline{University of California, San Diego}
\centerline{\tt pfitzsim@ucsd.edu}

\bigskip

In a recent note  \ref{Cu}, N.\ Curien pointed out that the Strong Law of Large Numbers is a direct consequence of the following fact.

\proclaim{Proposition}  If $(X_k)$ are iid integrable random variables with $\EEE[X_k]>0$, then
$$
\inf_n S_n>-\infty,\qquad\hbox{a.s.}
\leqno(1)
$$
\endproclaim

\nin Here, $S_n=X_1+\cdots+X_n$.  
In particular, (1) implies that
$$
\liminf_n {S_n\over n}\ge 0,\qquad\hbox{a.s.},
\leqno(2)
$$
which is the key to the SLLN.
[Indeed, for  $(X_k)$ an iid sequence with finite mean $\mu$, fix $c<\mu$ and apply (2) to $(X_k-c)$ to see that
$$
\liminf_n{S_n-nc\over n}\ge 0,\qquad \hbox{a.s.}
$$
 That is,
$$
\liminf_n {S_n\over n}\ge c,\qquad\hbox{a.s.}
$$
It follows that $\liminf_n S_n/n\ge\mu$. 
Apply the same reasoning to $(-X_k)$ to see that $\limsup_n S_n/n\le\mu$.]

The main point of \ref{C}, however,  is a quick proof of (1), using the notion of descending ladder times, time-reversal duality, and Wald's first identity.

Here we present a different elementary proof of the Proposition. Let $(X_k)$ be an iid sequence with finite mean $\mu=\EEE[X_k]>0$, and define $S_n$ as before. For positive integer $n$ define 
$$
J_n:=\min(S_1,S_2,\ldots,S_n),
$$
and note that each $J_n$ is integrable and $n\mapsto J_n$ is decreasing. Evidently,
$$
J_n=X_1+\min(0,S_2-S_1,S_3-S_1,\ldots,S_n-S_1)=X_1+\min(0,J_{n-1}'),\leqno(3)
$$
where $J_{n-1}'=\min(S_2-S_1,S_3-S_1,\ldots,S_n-S_1)$ has the same distribution as $J_{n-1}$. Therefore
$$
\EEE[J_n]=\mu+\EEE[\min(0,J_{n-1})], \qquad n\ge 2.
$$
Subtracting $\EEE[J_{n-1}]$ from both sides we find that 
$$
\EEE[J_n-J_{n-1}]=\mu-\EEE[J_{n-1}^+], \qquad n\ge 2,
$$
where $b^+$ denotes the positive part of $b$. Telescoping:
$$
\sum_{k=1}^{n}\EEE[J_k^+]=n \mu+\EEE[J_1-J_{n+1}]\ge n \mu ,\qquad n\ge 2.
$$
It follows that
$$
\liminf_n n^{-1}\sum_{k=1}^{n}\EEE[J_k^+]\ge\mu,
$$
and so, since $n\mapsto \EEE[J_n^+]$ is monotone decreasing, that 
$$
\EEE[J^+_n]\ge\mu,\qquad\forall n\ge 2.\leqno(4)
$$
Write $J_\infty:=\downarrow \lim_nJ_n=\inf_{k\ge 1}S_k$. By (4) and monotone convergence
$$
\EEE[J_\infty^+]\ge\mu>0,
$$
and in particular, $\PPP[J_\infty>-\infty]\ge\PPP[J_\infty>0]>0$. By Kolmogorov's Zero-One Law, $\PPP[J_\infty>-\infty]=1$. \qed
\bigskip

\nin{\bf Remark.} The reader will have noticed that the argument above uses the mutual independence of the $X_k$ only to make the Zero-One law available. In short, the conclusion of the Proposition  remains valid if $(X_k)$ is a stationary ergodic sequence of integrable random variables. The conclusion that $\lim_n \displaystyle{S_n\over n}=\mu$ a.s. (the Ergodic Theorem) then follows as before. 

\bigskip

\nin{\bf Added Remarks.} (a) Shortly after I posted this (10 November 2021), N.~Curien apprised me of a blog post \ref{D} containing an argument similar  to that presented here (and ascribed by R. Douc to Bernard Delyon)  showing that $\displaystyle{\liminf_n{S_n\over n}\ge 0}$ a.s. The ``stochastic recursion" approach exemplified by (3) is not new, dating back at least as far as A.~Garsia's proof \ref{G} of the Hopf maximal inequality, and used masterfully (in operator form) by J.~Neveu in his treatment of the Chacon-Ornstein ratio ergodic theorem. 
A wonderful discussion of such things can be found in \ref{S}.
\medskip

(b) C.W. Chin has noted in \ref{Ch} that even without ergodicity, the conclusion $\EEE[\inf_nS_n]\ge\mu$ for stationary  $(X_n)$ with $\mu=\EEE[X_1]>0$ yields the inequality $\displaystyle{\liminf_n{S_n\over n}}\ge\EEE[X_1|\II]$ a.s. ($\II$ denoting the invariant $\sigma$-field), and then the full Birkhoff ergodic theorem.
\bigskip

Here is some detail on this point. Replacing $X_k$ by $X_k-\EEE[X_1|\II]$ (still stationary) we can assume that $\EEE[X_k|\II]=0$ and show that $\lim_nS_n/n=0$ a.s. As usual, it is enough to show that $\liminf_nS_n/n\ge 0$, a.s. Fix $\epsilon>0$ and suppose that the event $A:=\{\liminf_nS_n/n<-\epsilon\}\in\II$ has strictly positive probability.  Then the stationary sequence defined by $\tilde X_k:=1_A\cdot(X_k+\epsilon/2)$ satisfies
$$
\EEE[\tilde X_k]=\EEE[X_k; A]+{\epsilon\over 2}\cdot \PPP[A]={\epsilon\over 2}\cdot \PPP[A]>0,
$$
and so, writing $\tilde J_\infty:=\inf_n\tilde S_n$,
$$
\PPP[\tilde J_\infty>0]>0.
$$
But clearly $\{\tilde J_\infty>0\}\subset A$. Meanwhile, on $\{\tilde J_\infty>0\}$ we have
$$
\liminf_n{S_n\over n}\ge -\epsilon/2.
$$
This yields a contradiction of the very definition of $A$, unless $\PPP[A]=0$.

\centerline{\bf References}
\medskip
\frenchspacing

\item{[Ch]}
Chin, Calvin Wooyoung: A mass transport proof of the ergodic theorem, (2022).\hfill\break{\tt https://arxiv.org/abs/2203.09687}.
\smallskip

\item{[Cu]}
Curien, Nicolas: Yet another proof of the strong law of large numbers, (2021).\hfill\break  {\tt https://arxiv.org/abs/2109.04315} .
\smallskip

\item{[D]}
Douc, Randal: A short proof of the Strong Law of Large Numbers...,\hfill\break  {\tt https://wiki.randaldouc.xyz/doku.php?id=world:lln}
\smallskip

\item{[G]}
Garsia, Adriano M.: A simple proof of E. Hopf's maximal ergodic theorem. {\it J. Math. Mech.} {\bf 14} (1965) 381--382.
\smallskip

\item{[N]}
Neveu, Jacques:   The filling scheme and the Chacon-Ornstein theorem. { \it  Israel J. Math.} {\bf  33} (1979)
 368--377.
 \smallskip

\item{[S]}
Steele, J. Michael:   Explaining a mysterious maximal inequality --- and a path to the law of large numbers. {\it Amer. Math. Monthly} {\bf 122} (2015) 490--494.
\smallskip

\end